\documentclass[12 pt]{amsart}
\usepackage{amsfonts}
\usepackage{amssymb}
\usepackage{amsmath}
\usepackage{graphicx}
\usepackage{color}

\theoremstyle{definition}
\newtheorem{defs}{Definition}[section]
\newtheorem{APdefs}{Definition}

\theoremstyle{plain}

\newtheorem{lem}[defs]{Lemma}
\newtheorem{ths}[defs]{Theorem}
\newtheorem{cor}[defs]{Corollary}
\newtheorem*{NNths}{Theorem}

\newtheorem{APprop}[APdefs]{Proposition}

\theoremstyle{remark}
\newtheorem{rem}[defs]{Remark}
\newtheorem*{NNrem}{Remark}
\newtheorem{Example}[defs]{Example}
\newtheorem{APrem}[APdefs]{Remark}

\def\cv{``$\circ$''}

\title[Configurations of saddle connections]
{Configurations of saddle connections of quadratic differentials on
$\mathbb{CP}^1$ and on 
hyperelliptic Riemann surfaces}
\author{Corentin Boissy}
\address{IRMAR,
Universit\'e Rennes 1,
Campus de Beaulieu,
35042 Rennes cedex, France}
\email{corentin.boissy@univ-rennes1.fr}


\begin{document}

\begin{abstract}
Configurations of rigid collections of saddle connections are connected component invariants for strata of the moduli space of quadratic differentials. They have been classified for strata of Abelian differentials by Eskin, Masur and Zorich. Similar work for strata of quadratic differentials has been done in Masur and Zorich, although in that case the connected components were not distinguished. 

We classify the configurations for quadratic
differentials on $\mathbb{CP}^1$ and on hyperelliptic connected components of
the moduli space of quadratic differentials. We show that, in genera greater than five, any configuration that appears in the hyperelliptic connected component of a stratum also appears in the non-hyperelliptic one. 
\end{abstract}

\maketitle

\setcounter{tocdepth}{1}
\tableofcontents

\section{Introduction}
We study flat surfaces having conical singularities of angle integer multiple of $\pi$ and $\mathbb{Z}/2\mathbb{Z}$ linear holonomy. The moduli space of such surfaces is isomorphic to the moduli space of quadratic differentials on
Riemann surfaces and is naturally stratified. Flat surfaces corresponding to squares of Abelian differentials are often called \emph{translation surfaces}. Flat surfaces appear in the study of billiards in rational polygons since these can be "unfolded" to give a translation surface (see \cite{KaZe}).

A sequence of quadratic differentials or Abelian differentials leaves any compact set of a stratum when the length of a saddle connection tends to zero. This might force some other saddle connections to shrink. In the case of an Abelian differential this correspond to homologous saddle connections. In the general case of quadratic differentials, the corresponding collections of saddle connections on a flat surface are said to be \emph{\^homo\-logous}\footnote{The corresponding cycles are in fact homologous on the canonical double cover of $S$, usually denoted as $\widehat{S}$, see section \ref{homologous}.}
 (pronounced ``hat-homologous''). 
According to Masur and Smillie \cite{MS} (see also \cite{EMZ,MZ}), a ``typical degeneration'' corresponds to the case when all the ``short'' saddle connections are pairwise \^homo\-logous). Therefore the study of  configurations of \^homo\-logous saddle connections (or homologous saddle connection in the case of Abelian differential) is related to the study of the compactification of a given stratum. A  configuration of \^homo\-logous saddle connections on a generic surface is also a natural invariant of a connected component of the ambient stratum.

In a recent article, Eskin, Masur and Zorich \cite{EMZ} study collections of homologous saddle connections for Abelian differentials. They describe configurations for each connected component of the strata of Abelian differentials. Collections of \^homo\-logous saddle connections are studied for quadratic differentials by  Masur and Zorich \cite{MZ}:  they describe all the configurations that can arise in any given stratum of quadratic differentials, but they do not distinguish connected components of such strata.

According to Lanneau \cite{La2}, the non-connected strata of quadratic differentials admit exactly two  connected components. They are of one of the following two types:
\begin{itemize}
\item ``hyperelliptic'' stratum: the stratum admits a connected component that consists of hyperelliptic quadratic differentials (note that some of these strata are connected).
\item exceptional stratum: there exist four non-connected strata that do not belong to the previous case.
\end{itemize}

In this article, we give the classification of the configurations that appear in the hyperelliptic connected components (Theorem \ref{confighyp}). 
This gives therefore a necessary condition for a surface to be in a hyperelliptic connected component.  Unfortunately, this is not a sufficient condition since, as we show, any configuration that appears in a hyperelliptic connected component  also appears in the other component of the stratum when the  genus is greater than $5$ (Theorem \ref{config:non:hyp}). We address the description of configurations for exceptional strata to a next article.

We deduce configurations for hyperelliptic components from configurations for strata of quadratic differentials on $\mathbb{CP}^1$ (Theorem \ref{prop:config}). Configurations for $\mathbb{CP}^1$ are deduced from general results on configurations that appear in~\cite{MZ}. 
Note that these configurations are needed in the study of asymptotics in billiards in polygons with ``right'' angles \cite{AEZ}. For such a polygon, there is a simple unfolding procedure that consists in gluing along their boundaries two copies of the polygon. This gives a flat surface of genus zero with conical singularities, whose angles are multiples of $\pi$ (\emph{i.e.} a quadratic differential on $\mathbb{CP}^1$). Then a generalized diagonal or a periodic trajectory in the polygon gives a saddle connection on the corresponding flat surface.

We also give in appendix an explicit formula that gives a relation between the genus of a surface and the ribbon graph of connected components associated to a collection of \^homo\-logous saddle connections.

Some particular splittings are sometimes used to compute the closure of $\textrm{SL}(2,\mathbb{R})$-orbits of surfaces (see \cite{McMull, HLM}). These splittings of surfaces can be reformulated as configurations of homologous or \^homo\-logous saddle connections on these surfaces. It would be interesting to find some configurations that appear in \emph{any} surface of a connected component of a stratum, as was done in \cite{McMull}.

\subsubsection*{Acknowledgements} I would like to thank Anton Zorich for encouraging me to write this paper, and for many discussions.
I also thank Erwan Lanneau, for his many comments.

\subsection{Basic definitions}\label{bas_def}
Here we first review standart facts about moduli spaces of quadratic differentials. We refer to \cite{Hubbard:Masur,Masur82,Veech82} for proofs and details, and to \cite{MT,Z} for general surveys.

Let  $S$  be a  compact  Riemann surface  of  genus  $g$. A  quadratic differential $q$  on $S$ is  locally given by  $q(z)=\phi(z)dz^2$, for $(U,z)$ a local chart with  $\phi$ a meromorphic function with at most simple poles. 
We define the poles and zeroes of $q$ in a local chart to be the poles and zeroes of the corresponding meromorphic function $\phi$. It is easy to check that they do not depend on the choice of the local chart. Slightly abusing  vocabulary, a pole will be referred to as a zero of order $-1$, and a marked point will be referred to as a zero of order $0$. An Abelian differential on $S$ is a holomorphic 1-form.

Outside its poles and zeroes,  $q$ is locally the square of an Abelian differential. Integrating this 1-form gives a natural atlas such that  the transition  functions are  of the kind $z\mapsto \pm z+c$. Thus $S$  inherits a flat metric with singularities, where  a zero  of order  $k\geq -1$ becomes  a conical singularity of angle $(k+2)\pi$.  The flat metric has trivial holonomy if and only if $q$ is globally the square of an Abelian differential. If not, then the holonomy  is   $\mathbb{Z}/2\mathbb{Z}$  and  $(S,q)$  is sometimes called  a \emph{half-translation} surface since the transitions functions are either translations, either half-turns. In order to simplify the notation, we will usually denote by $S$ a surface with a flat structure.

We associate to a quadratic differential the set $\{k_1,\ldots,k_r\}$ of orders of its  poles and zeros.  The Gauss-Bonnet  formula asserts that $\sum_i  k_i=4g-4$. Conversely, if we fix a collection $\{k_1,\dots,k_r\}$ of integers greater than or equal to $-1$ satisfying the previous equality, we denote by $\mathcal{Q}(k_1,\ldots,k_r)$ the (possibly empty) moduli space of quadratic differentials which are not globally the square  of any Abelian differential, and having $\{k_1,\ldots,k_r\}$ as orders of poles and zeroes .  It is well known that $\mathcal{Q}(k_1,\ldots,k_r)$ is a complex analytic orbifold, which  is usually called a \emph{stratum} of the moduli space of quadratic differentials. We mostly restrict ourselves to the subspace $\mathcal{Q}_1(k_1,\ldots,k_r)$  of area one surfaces, where the area is given by the flat metric.  In a similar way, we denote by $\mathcal{H}_1(n_1,\dots,n_s)$ the moduli space of Abelian differentials  of area $1$ 
having zeroes of degree $\{n_1,\ldots,n_s\}$, where $n_i\geq 0$ and $\sum_{i=1}^s n_i=2g-2$.

A saddle  connection is  a geodesic segment (or geodesic loop) joining two singularities (or a singularity to itself) with no  singularities in its interior. Even if $q$ is not  globally a square of an Abelian differential we can find a square   root   of  it   along  the   saddle connection.  Integrating  it along  the  saddle  connection  we get  a complex number  (defined up  to multiplication by  $-1$). Considered as a planar vector, this complex number represents the affine holonomy vector along the saddle connection. In particular, its euclidean length is the modulus of its  holonomy vector. Note that  a saddle connection persists under  small deformation of the  surface. 

Local coordinates on a stratum of Abelian differentials are obtained by integrating the holomorphic 1-form along a basis of the relative homology $H_1(S,{sing},\mathbb{Z})$, where \emph{sing} is the set of conical singularities. Equivalently, this means that local coordinates are defined by the relative cohomology $H^1(S,{sing},\mathbb{C})$.

Local coordinates in a stratum of quadratic differentials are obtained by the following way: one can naturally associate to a quadratic differential $(S,q)\in
\mathcal{Q}(k_1,\dots,k_r)$ a double cover $p:\widehat{S}\rightarrow S$ such that $p^* q$ is the square of an Abelian differential $\omega$.  The surface $\widehat{S}$ admits a natural involution $\tau$, that induces on the relative cohomology $H^1(S,{sing},\mathbb{C})$ an involution $\tau^*$. It decomposes $H^1(S,{sing},\mathbb{C})$ into an invariant subspace $H^{1}_+(S,{sing},\mathbb{C})$ and an anti-invariant subspace $H^1_-(S,{sing},\mathbb{C})$. One can show that the anti-invariant subspace $H^1_-(S,{sing},\mathbb{C})$ gives local coordinates for the stratum $\mathcal{Q}(k_1,\ldots,k_r)$. It is well known that Lebesgue measure on these coordinates defines a finite measure $\mu$ on the stratum $\mathcal{Q}_1(k_1,\ldots,k_r)$.

A hyperelliptic quadratic differential is a quadratic differential such that there exists an orientation preserving involution $\tau$ with $\tau^*q=q$  and such that  $S/\tau$ is a  sphere.  We can construct families of hyperelliptic quadratic differentials by the following way: to all quadratic differentials on $\mathbb{CP}^1$, we associate a double covering ramified over some singularities satisfying some fixed combinatorial conditions. The resulting Riemann surfaces  naturally carry  hyperelliptic quadratic differentials.
 
Some strata  admit an entire connected   component  that  is made of   
hyperelliptic  quadratic differentials. These components arise from the previous
construction and have been classified by M.~Kontsevich and A.~Zorich in case of  Abelian differentials \cite{KoZo} and by E. Lanneau in case of quadratic differentials \cite{La1}. 

\begin{NNths}[M. Kontsevich,  A. Zorich]
The strata of Abelian differentials having a hyperelliptic connected component are the following ones.
\begin{enumerate}
\item $\mathcal{H}(2g-2)$, where $g\geq 1$. It arises from
	$\mathcal{Q}(2g-3,-1^{2g+1})$. The ramifications points are located 
over all the singularities.
\item $\mathcal{H}(g-1,g-1)$, where $g\geq 1$. It arises from
	$\mathcal{Q}(2g-2,-1^{2g+2})$. The ramifications points are located
over all the poles.
\end{enumerate}
In the above presented list, the strata $\mathcal{H}(0)$, $\mathcal{H}(0,0)$,  $\mathcal{H}(1,1)$ and $\mathcal{H}(2)$ are the ones that are connected.
\end{NNths}

\begin{NNths}[E. Lanneau]
The  strata of quadratic differentials that have a hyperelliptic  connected
component  are the following ones.
\begin{enumerate}
\item  $\mathcal{Q}(2(g-k)-3,2(g-k)-3,2k+1,2k+1)$ where $k\geq -1$,  $g\geq 1$ and $g-k\geq 2$.  It arises from $\mathcal{Q}(2(g-k)-3,2k+1,-1^{2g+2})$. The ramifications points are located over $2g+2$ poles.  
\item  $\mathcal{Q}(2(g-k)-3,2(g-k)-3,4k+2)$ where  $k\geq 0$, $g\geq 1$ and $g-k\geq 1$. It arises from $\mathcal{Q}(2(g-k)-3,2k,-1^{2g+1})$. The ramifications points are located over $2g+1$ poles and over the zero of order $2k$.
\item $\mathcal{Q}(4(g-k)-6,4k+2)$ where  $k\geq 0$, $g\geq 2$ and $g-k\geq 2$.  It arises from $\mathcal{Q}(2(g-k)-4,2k,-1^{2g})$. The ramifications points are located over all the singularities 
\end{enumerate}
In the above presented list, the strata $\mathcal{Q}(-1,-1,1,1)$, $\mathcal{Q}(-1,-1,2)$, $\mathcal{Q}(1,1,1,1)$, $\mathcal{Q}(1,1,2)$ and $\mathcal{Q}(2,2)$ are the ones that are connected.
\end{NNths}

\subsection{\^Homo\-logous saddle connections}\label{homologous}

Let $S\in \mathcal{Q}(k_1,\dots,k_r)$ be a flat surface and let us denote by $p:\widehat{S}\rightarrow S$ its canonical double cover and by $\tau$ the corresponding involution. Let $\Sigma$ denote the set of singularities of $S$ and let $\widehat{\Sigma}=p^{-1}(\Sigma)$.

To an oriented saddle connection $\gamma$ on $S$, one can associate $\gamma_1$ and $\gamma_2$ its preimages by $p$. If the relative cycle $[\gamma_1]$ satisfies $[\gamma_1]=-[\gamma_2] \in H_1(\widehat{S},\widehat{\Sigma},\mathbb{Z})$, then we define $[\tilde{\gamma}]=[\gamma_1]$. Otherwise, we define  $[\tilde{\gamma}]=[\gamma_1]-[\gamma_2]$. Note that in all cases, the cycle $[\tilde{\gamma}]$ is anti-invariant with respect to the involution~$\tau$.

\begin{defs}\label{def:homo}
Two saddle connections $\gamma$ and $\gamma^{\prime}$ are \^homo\-logous if $[\tilde{\gamma}]=\pm[\tilde{\gamma^{\prime}}]$.
\end{defs}

\begin{figure}[htb]
\begin{center}
\input{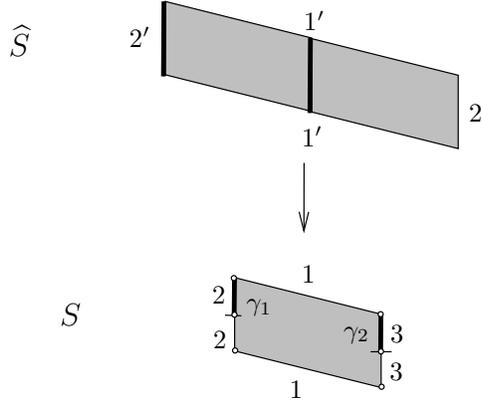}
\caption{An unfolded flat surface $S$ with two \^homo\-logous saddle connections $\gamma_1$ and $\gamma_2$.}
\label{exemple_trivial}
\end{center}
\end{figure}

\begin{Example}
Consider the flat surface $S\in \mathcal{Q}(-1,-1,-1,-1)$ given in Figure \ref{exemple_trivial} (a ``pillowcase''), it is easy to check from the definition that $\gamma_1$ and $\gamma_2$ are \^homo\-logous since the corresponding cycles for the double cover $\widehat{S}$ are homologous.
\end{Example}

\begin{Example}
Consider the flat surface given in Figure \ref{exempleplat} (at the end of section \ref{homologous}), the reader can check that the saddle connections $\gamma_1$, $\gamma_2$ and $\gamma_3$ are pairwise \^homo\-logous.
\end{Example}

\begin{NNths}[H. Masur, A. Zorich]
Consider two distinct saddle connections $\gamma,\gamma^{\prime}$ on a half-translation surface. The following assertions are equivalent:
\begin{itemize}
\item
The two saddle connections $\gamma$ and $\gamma^{\prime}$ are \^homo\-logous.
\item
The ratio of their lengths is constant under any small deformation of the surface inside the ambient stratum.
\item
They  have no interior intersection and one of the
connected  components  of $S\backslash\{\gamma\cup \gamma^{\prime}\}$
has  trivial  linear holonomy.   
\end{itemize}
Furthermore, if  $\gamma$ and $\gamma^{\prime}$ are \^homo\-logous, then the ratio of their lengths belongs to $\{\frac{1}{2},1,2\}$.
\end{NNths}

Consider a set of  \^homo\-logous  saddle connections $\gamma=\{\gamma_1,\dots,\gamma_s\}$  on  a flat surface $S$. Slightly abusing notation, we will denote by $S\backslash   \gamma$ the subset $S\backslash\bigl( \cup_{i=1}^s \gamma_i \bigr)$. This subset is a finite union of connected
half-translation surfaces with boundaries. 

\begin{defs}
Let $S$ be a flat surface and $\gamma=\{\gamma_1,\ldots,\gamma_s\}$ a collection of \^homo\-logous saddle connections. 
The \emph{graph of connected components}, denoted by $\Gamma(S,\gamma)$, is the graph defined by the following way:
\begin{itemize}
\item The vertices are  the  connected components   of  $S\backslash   \gamma$,  labelled   by  ``$\circ$''   if  the corresponding surface is a cylinder,  by ``$+$'' if it has trivial holonomy (but is  not a cylinder), and  otherwise by ``$-$'' if it  has non-trivial holonomy.
\item The  edges are  given  by  the saddle  connections  in $\gamma$. Each
$\gamma_i$ is located on the  boundary of one  or two connected  components of $S\backslash  \gamma$. In  the  first case  it  becomes an  edge joining  the
corresponding vertex to itself. In the second case, it becomes an edge
joining the two corresponding vertices. 
\end{itemize}
\end{defs}

In  \cite{MZ}, Masur and Zorich describe the set of all possible graphs of connected components for a quadratic differential.
This set is roughly given by Figure~\ref{fig:classification:of:graphs}, where dot lines are chains of ``$+$'' and   ``$\circ$''  vertices of valence two. The next theorem gives a more precise statement of this description. It can be skipped in a first reading.

\begin{NNths}[H. Masur, A. Zorich]
Let  $(S,q)$ be quadratic  differential;  let
$\gamma$ be  a  collection  of  \^homo\-logous  saddle  connections
$\{\gamma_1, \dots,  \gamma_n\}$,  and  let $\Gamma(S,\gamma)$ be
the graph  of  connected  components  encoding  the decomposition
$S\setminus (\gamma_1\cup\dots\cup\gamma_n)$.

The graph $\Gamma(S,\gamma)$ either has one of the basic types listed
below or can  be  obtained  from one of these graphs  by  placing
additional  ``$\circ$''-vertices  of valence two at any subcollection  of
edges  subject  to  the  following  restrictions.   At  most  one
\cv-vertex may be placed at the same edge;  a \cv-vertex cannot
be placed at an  edge adjacent to a \cv-vertex of valence  $3$ if
this is the edge separating the graph.

The     graphs     of     basic      types,      presented     in
Figure~\ref{fig:classification:of:graphs},  are  given   by   the
following list:
\begin{itemize}
\item[a)]
An  arbitrary  (possibly  empty)  chain  of  ``$+$''-vertices  of
valence two bounded by a pair of ``$-$''-vertices of valence one;
\item[b)]
A  single loop  of  vertices of valence  two  having exactly  one
``$-$''-vertex and arbitrary number of ``$+$''-vertices (possibly
no ``$+$''-vertices at all);
\item[c)]
A single  chain and  a single loop joined at  a vertex of valence
three. The graph has  exactly  one ``$-$''-vertex of valence one;
it is  located at the  end of  the chain. The  vertex of  valence
three is either a ``$+$''-vertex, or a \cv-vertex  (vertex of the
cylinder  type).  Both the  chain,  and  the  cycle  may  have in
addition an arbitrary  number  of  ``$+$''-vertices  of  valence two
(possibly no ``$+$''-vertices at all);
\item[d)]
Two nonintersecting  cycles joined by a  chain. The graph  has no
``$-$''-vertices. Each  of the two  cycles has a single vertex of
valence three (the one where the chain is attached to the cycle);
this vertex  is either a  ``$+$''-vertex or a \cv-vertex. If both
vertices of valence three are \cv-vertices, the chain joining two
cycles  is  nonempty:   it   has  at  least  one  ``$+$''-vertex.
Otherwise, each  of the cycles and  the chain may  have arbitrary
number  of   ``$+$''-vertices   of   valence   two  (possibly  no
``$+$''-vertices of valence two at all);
\item[e)]
``Figure-eight'' graph: two  cycles joined at a vertex of valence
four, which  is either a  ``$+$''-vertex or a \cv-vertex. All the
other vertices (if  any) are the ``$+$''-vertices of valence two.
Each  of  the  two  cycles  may  have arbitrary  number  of  such
``$+$''-vertices of  valence two (possibly no ``$+$''-vertices of
valence two at all).
\end{itemize}

Each graph listed above corresponds to some flat surface $S$ and
to some collection of saddle connections $\gamma$.
   \begin{figure}[htb]
%
\includegraphics{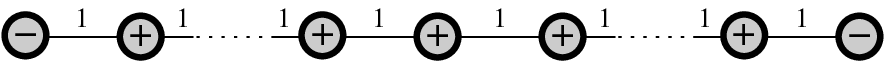}
   %
\includegraphics{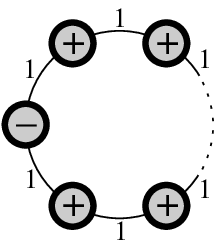}
   %
\includegraphics{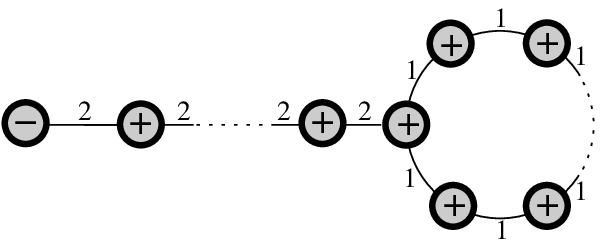}
\includegraphics{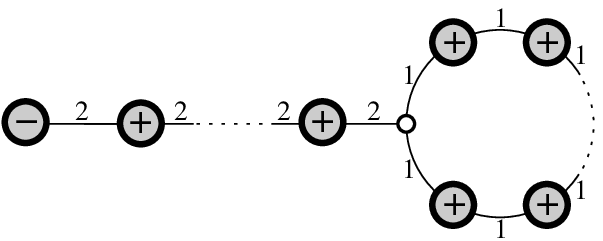}
   %
\includegraphics{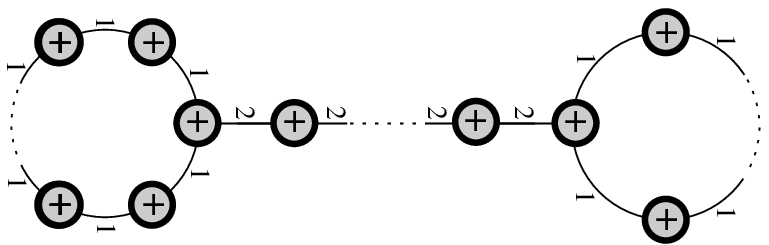}
\includegraphics{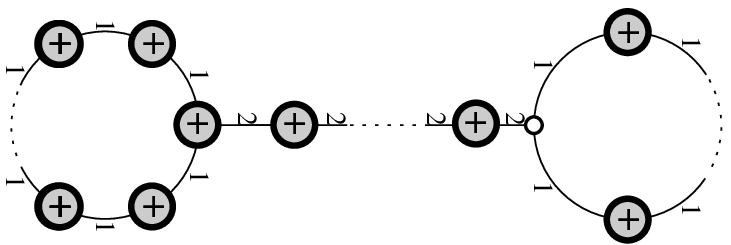}
\includegraphics{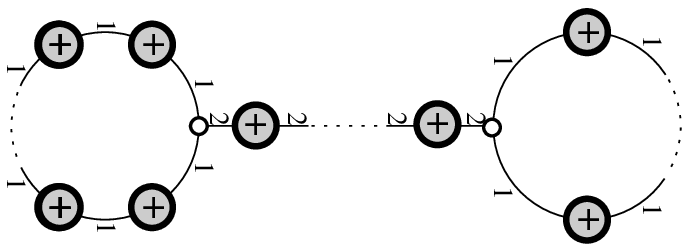}
   %
\includegraphics{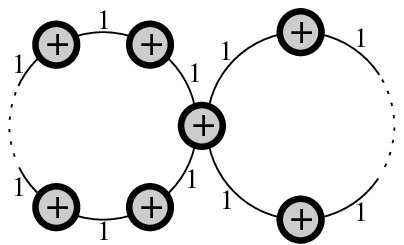}
\includegraphics{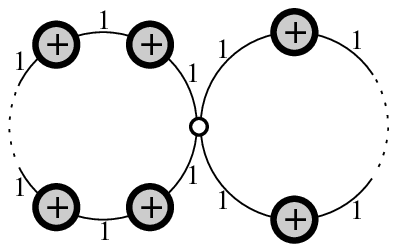}
%
%
%
\vspace{490bp} 
\begin{picture}(0,0)(0,0)
\put(-140,470.5)
{
\begin{picture}(0,0)(0,0)
\put(38,-29){a)}
\put(285,-53){b)}
\put(132,-130){c)}
\put(132,-362){d)}
\put(132,-450){e)}
\end{picture}}
\end{picture}
\caption{
\label{fig:classification:of:graphs}
Classification of admissible graphs.
}
\end{figure}

\end{NNths}

\begin{rem}
Two \^homo\-logous saddle connections are not necessary of the same length.
The additional parameters $1$ or $2$ written along the vertices in Figure \ref{fig:classification:of:graphs} represent the lengths of the saddle connections in the collection $\gamma=\{\gamma_1,\ldots,\gamma_s\}$ after suitably rescaling the surface.
\end{rem}

Each connected  component of $S\backslash \gamma$ is a non-compact
surface which can be naturally compactified (for example considering
the distance induced by the flat metric on a connected component of
$S\backslash \gamma$, and the corresponding completion).
We denote this compactification by $S_j$. We warn the reader that
$S_j$ might  differ from the  closure of
the component in the surface $S$: for example, if $\gamma_i$ is on the
boundary of just one connected  component $S_j$ of $S\backslash\gamma$, then
the compactification of $S_j$ carries  two copies of $\gamma_i$ in its
boundary,  while  in  the  closure  of  the corresponding connected component of $S\backslash \gamma$,   these  two  copies  are
identified.  The  boundary of each  $S_i$ is a union of saddle connections;
it has one or  several connected components.
Each of them is  homeomorphic to $\mathbb{S}^1$ and therefore defines a
cyclic order in the set of boundary saddle connections.
Each consecutive pair of saddle connections for that cyclic order
defines a \emph{boundary singularity} with an associated angle which
is a integer multiple of $\pi$ (since the boundary saddle
connections are parallel). The surface with boundary $S_i$ might have
singularities in its interior. We call them \emph{interior
singularities}.

\begin{defs}
Let $\gamma=\{\gamma_1,\dots,\gamma_r\}$ be a maximal collection of
\^homo\-logous saddle connections. Then a \emph{configuration} is the
following combinatorial data:
\begin{itemize}
\item The graph $\Gamma(S,\gamma)$.
\item For each vertex of this graph, a permutation of the edges
adjacent to the vertex (encoding the cyclic order of the
saddle connections
on each connected component of the boundary of $S_i$).
\item For each pair of consecutive elements in that cyclic
order, an integer $k\geq 0$ such that the angle between the two corresponding
saddle connections is $(k+1)\pi$. This integer will be referred as the \emph{order of the boundary singularity}. 
\item For each $S_i$, a collection of integers corresponding to the orders of
the interior singularities of $S_i$.
\end{itemize}
\end{defs}

Following \cite{MZ}, we will encode the permutation of the edges adjacent to each vertex by a \emph{ribbon graph}.

\begin{figure}[htpb]
\begin{center}
\input{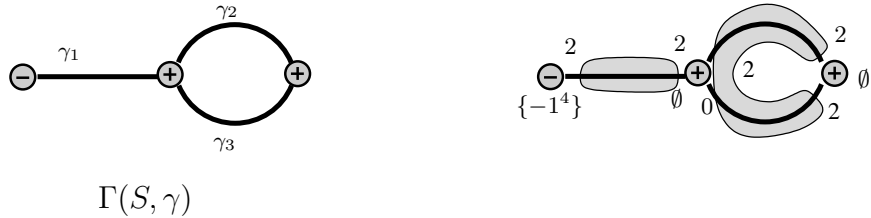}
\caption{An example of configuration.}
\label{exemple}
\end{center}
\end{figure}

\begin{Example}
Figure \ref{exemple} represents a configuration on a flat surface. The corresponding collection $\{\gamma_1,\gamma_2,\gamma_3\}$ of \^homo\-logous saddle connections decomposes the surface into three connected components.  The first connected component  \includegraphics{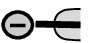} has four interior singularities of order $-1$, and its boundary consists of a single saddle connection with the corresponding boundary singularity of angle $(2+1)\pi=3\pi$. The second connected component \includegraphics{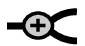}  has no interior singularities. It has two boundary components, one consisting of a single saddle connection with corresponding singularity of angle $(2+1)\pi$, and the other consists of a union of two saddle connections with corresponding boundary singularities of angle $(0+1)\pi$ and $(2+1)\pi$. The last connected component~\includegraphics{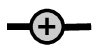} has no interior singularities, and admits two boundary components that consists each of a  single saddle connection with corresponding  boundary singularities of angles $(2+1)\pi$.

Figure \ref{exempleplat} represents a flat surface with a collection of three \^homo\-logous saddle connections realizing  this configuration.
\end{Example}

\begin{figure}[htpb]
\begin{center}
\input{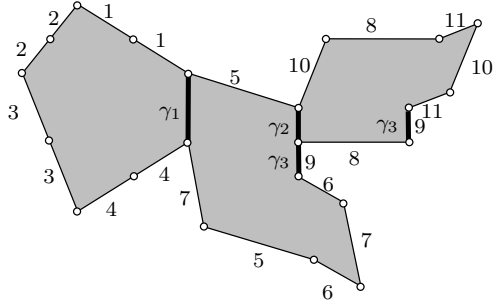}
\caption{Unfolded flat surface realizing configuration of Figure \ref{exemple}.}
\label{exempleplat}
\end{center}
\end{figure}

\begin{rem}\label{generique}
When describing the configuration of a collection of \^homo\-logous saddle connections $\gamma=\{\gamma_1,\ldots,\gamma_r\}$, we will always assume that the quadratic differential is \emph{generic}, and therefore, each saddle connection parallel to the $\gamma_i$ is actually \^homo\-logous to the $\gamma_i$ (see \cite{MZ}).
\end{rem}

\begin{rem}
A maximal collection of \^homo\-logous saddle connections and the associated configuration persist under any small deformation of the flat surface inside the ambient stratum. They also persist under the well know $SL(2,\mathbb{R})$ action on the stratum which is ergodic with respect to the Lebesgue measure $\mu$  (see \cite{Masur82,Veech82,Veech86}). 
Hence, all admissible configurations that exists in a connected component are realized in a generic surface of that component. Furthermore, the number of collections realizing a given configuration in a generic surface has quadratic asymptotics (see \cite{EM}).
\end{rem}

\section{Configurations for the Riemann sphere}
In this section we describe all  admissible configurations of \^homo\-logous
saddle connections that arise on $\mathbb{CP}^1$.
To avoid confusion of notation, we specify the following convention: we denote by $\{k_1^{\alpha_1},\ldots,k_r^{\alpha_r}\}$ the set  with multiplicity $\{k_1,k_1, \ldots,k_r\}$, where $\alpha_i$ is the multiplicity of $k_i$.
We assume that $k_i \neq k_j$ for $i\neq j$. For example the notation $\mathcal{Q}(1^2,-1^6)$ means $\mathcal{Q}(1,1,-1,-1,-1,-1,-1,-1)$.

Let $\mathcal{Q}(k_1^{\alpha_1},\ldots,k_r^{\alpha_r}, -1^s)$ be a stratum of quadratic differentials on $\mathbb{CP}^1$ different from $\mathcal{Q}(-1^4)$. 
We give in the next example four families of admissible configurations for this stratum. In the next example, $\gamma$ is always assumed to be a  maximal collection of \^homo\-logous saddle connections. We give in Table \ref{graphcp1} the corresponding graphs and ``topological pictures''. The existence of each of these configurations is a direct consequence of the  Main Theorem of \cite{MZ}.

\begin{Example}\label{exCP1}

\begin{itemize}
\item[a)]
Let $\{k,k^{\prime}\}\subset \{k_1^{\alpha_1},\ldots,k_r^{\alpha_r}, -1^s\}$ be  an unordered pair of integers with $(k,k^{\prime})\neq (-1,-1)$. The set $\gamma$ consists of a single saddle connection joining a singularity of order
$k$ to a distinct singularity of order $k^{\prime}$.
\item[b)]
Let $\{a_1,a_2\}$ be an unordered pair of strictly positive integers such that $a_1+a_2=k\in\{k_1, \ldots,k_r\}$ (with $k\neq 1$), and let $A_1 \sqcup A_2$ be a partition of $\{k_1^{\alpha_1},\ldots,k_r^{\alpha_r}\} \backslash \{k\}$. The set $\gamma$ consists of a simple saddle connection that decomposes the sphere into two one-holed spheres $S_1$ and $S_2$, such that each $S_i$ has interior singularities of positive order given by $A_i$ and $s_i=(\sum_{a\in A_i}a)+a_i+2$ poles, and has a~single boundary singularity of order $a_i$.
\item[c)]
Let $\{a_1,a_2\}\subset \{k_1^{\alpha_1},\ldots,k_r^{\alpha_r}\}$ be an unordered  pair of integers. Let $A_1 \sqcup A_2$ be a partition of $\{k_1^{\alpha_1},\ldots,k_r^{\alpha_r}\}\backslash \{a_1,a_2\}$. The set $\gamma$ consists of two closed saddle connections that decompose the sphere into two one-holed spheres $S_1$ and $S_2$ and a  cylinder,  and such that each $S_i$ has interior singularities of positive orders given by $A_i$ and  $s_i=(\sum_{a\in A_i}a)+a_i+2$ poles and has a boundary singularity of order $a_i$.
\item[d)]
Let $k\in \{k_1,\ldots,k_r\}$.  The set $\gamma$ is a pair of saddle connections of different lengths, and such that the largest one  starts and ends from a singularity of order $k$ and  decompose the surface into a one-holed sphere and a ``half-pillowcase'', while the shortest one joins a pair of poles and lies on the other end of the half-pillowcase.
\end{itemize}
\end{Example}

\begin{table}[htbp]
   \begin{center}
     \input{configcp1bis.pstex_t}
   \end{center}
   \caption{Configurations in genus zero}
   \label{graphcp1}
\end{table}

Now we will prove that the configuration described previously are the only ones for the stratum $\mathcal{Q}(k_1^{\alpha_1},\ldots,k_r^{\alpha_r}, -1^s)$.
We first start with several preliminary lemmas which are applicable to flat surfaces of arbitrary genus. Let $S$ be a generic flat surface of genus $g\geq 0$ in some stratum of quadratic differentials, and let $\gamma$ be a maximal collection of \^homo\-logous saddle connections on it. Taking the natural compactification of each connected  component of $S\backslash \gamma$, we get a collection $\{S_i\}_{i\in I}$ of compact surfaces with boundaries. The boundary of each $S_i$ is topologically a union of disjoint circles. We can glue a disc  to  each connected component of the boundary of $S_i$ and get a closed surface $\overline{S}_i$;  we denote by $g_i$ the genus of $\overline{S}_i$. 

\begin{lem}\label{genre}
        Let $g$ be the genus of $S$, then $g \geq \sum_{i\in I} g_i$.
\end{lem}
\begin{proof}
        For   each    $S_i$,   we   consider   a    collection   of
                  paths
        $(c_{i,1},\ldots,c_{i,2g_i})$   of  $S_i$  that   represent  a
        symplectic basis of $H_1(\overline{S}_i,\mathbb{R})$ and that avoid
        the boundary of  $S_i$. When we glue the  $\{S_i\}$ together, the $c_{i,j}$
        provides a  collection of cycles of $H_1(S,\mathbb{R})$.  It forms a
        symplectic family because two paths arising from two different
        surfaces do not  intersect each other. Therefore we get a  free family of
        $H_1(S,\mathbb{R})$, thus:
        $$2g=\dim \bigl(H_1(S,\mathbb{R})\bigr) \geq     \sum_{i\in I}
        \dim\bigl(H_1(\overline{S}_i,\mathbb{R})\bigr)=\sum_{i\in I} 2g_i.$$
\end{proof}

\begin{NNrem}
In the appendix, we will improve Lemma \ref{genre} and give an exact formula in terms of the graph $\Gamma(S,\gamma)$ and the ribbon graph.
\end{NNrem}

\begin{lem}\label{genre_plus}
        If  $S_{i_0}$ is  not a  cylinder and  has trivial  holonomy, then
        $g_{i_0}>0$.
\end{lem}
\begin{proof}
We recall that the initial collection of \^homo\-logous saddle connections is assumed to be maximal, therefore there is no  interior  saddle
connections \^homo\-logous   to  any   boundary   saddle  connection.    Let 
$\{k_1,\ldots,k_s\}$ be  the  orders  of  the  interior  conical singularities  of $S_{i_0}$  and  $\{l_1,\ldots,l_{s^{\prime}}\}$ be the orders of the boundary singularities.  Let  $X$ be the closed flat surface obtained by gluing $S_{i_0}$ and a copy of itself taken with opposite orientation along their boundaries. If $m$ denotes the number of connected components of  the  boundary of  $S_{i_0}$  and $g_X$ denotes the genus of $X$, one can see that $g_X=2g_{i_0}+m-1$.
 The singularities  of $X$ are of orders  $\{k_1,\ldots,k_s, k_1,\ldots,k_s,2l_1,\ldots,2l_{s^{\prime}}\}$.  Furthermore, $k_i,l_j$ are nonnegative integers since $X$ has trivial holonomy. Applying  the Gauss-Bonnet formula  for quadratic differentials,  one gets:
       \[g_X=1+\sum_{j=1}^s \frac{k_j}{2}
        +\sum_{i=1}^{s^{\prime}}   \frac{l_i}{2} 
        =2g_{i_0}+m-1\] 
which obviously gives 
        \[2g_{i_0}\geq 2-m + \sum_{i=1}^{s^{\prime}}    \frac{l_i}{2} .\]

To conclude, we need few elementary remarks (which are already written  in  \cite{MZ})  about  the order  of  the  conical singularities of the boundary:
\begin{itemize}
 \item[a)] 
If a connected  component of the boundary is just  a  single  saddle connection, then  the corresponding  angle cannot  be  $\pi$ otherwise the saddle connection would then be a boundary component of a cylinder.  Then the other boundary component of that cylinder would be a  saddle  connection \^homo\-logous to the previous one (see remark \ref{generique}). So $S_i$ would be that cylinder contradicting the hypothesis. 
 Furthermore, the holonomy of a path homotopic to the  saddle connection is trivial if  and  only  if  the conical  angle  of  the boundary  singularity is  an  odd multiple  of $\pi$ .  
 
Therefore that  angle  is greater  or  equal to $3\pi$, and hence, the corresponding order $l_j$ of the boundary singularity has order at least $2$.
\item[b)]
If a connected  component of the boundary is given  by  two  saddle  connections,  then  as before,  the  two  corresponding  conical angles cannot be  both equal to $\pi$ (otherwise $S_i$  would be a  cylinder)  and  are of  the same parity (otherwise $S_i$ would have nontrivial holonomy).
\end{itemize}

Now we complete the proof of the lemma. We recall that the vertex corresponding to $S_{i_0}$  in $\Gamma(S,\gamma)$ is of valence at most four, and hence $m\leq 4$. 
The case $m=1$ is trivial. If $m=2$ then there  is a connected component of  the boundary of $S_{i_0}$ with one or two saddle connections.  In both cases, the  remarks a) and b) imply  that  
$S_{i_0}$ admits a boundary singularity of order $l_1>0$, and therefore $2g_{i_0}\geq l_1>0$.

If  $m\in\{3,4\}$, then there are at least two boundary components that consist of a single saddle connection. From remark a), this implies that $S_{i_0}$ admits two boundary singularities $l_1$ and $l_2$ of order greater than or equal to two. Applying remarks a) and b) on the other boundary components, 
we show that $S_{i_0}$ admits at least an other boundary singularity of order $l_3>0$. Therefore 
\[2g_{i_0}>2-m+ l_1/2+l_2/2 \geq 4-m \geq 0.\] 
Finally, $g_{i_0} > 0$ and the lemma is proven.
\end{proof}

Now, we describe all the possible configurations when the genus $g$ of the surface is zero.

\begin{ths}\label{prop:config}
Let $\mathcal{Q}(k_1^{\alpha_1},\ldots,k_r^{\alpha_r}, -1^s)$ be a stratum of quadratic differentials on $\mathbb{CP}^1$ different from $\mathcal{Q}(-1^4)$, and let $\gamma$ be a maximal collection of \^homo\-logous saddle connections on a flat surface in this stratum. Then all possible configurations for $\gamma$ are the ones described in Example~\ref{exCP1}.
\end{ths}

\begin{proof}

It follows from Lemmas \ref{genre} and \ref{genre_plus} that  $\Gamma(S,\gamma)$ has no ``$+$'' component.  Furthermore, a loop of  the graph $\Gamma(S,\gamma)$ cannot  have   any  cylinder  since this   would  add  a  handle  to  the surface. Now using the description  from \cite{MZ} of admissible graphs (see Figure \ref{fig:classification:of:graphs}), we can list all  possible graphs. For each graphs, we now describe the corresponding admissible configurations.

a) A single  ``$-$'' vertex  of  valence two  and an  edge joining  it  to  itself.  \\ 
 This  can  represent  two possible  cases: either  the boundary  of the closure  of $S\backslash  \gamma$ has two connected components, or it has only one. In the  first case	each  connected component  of the boundary  is a single  saddle  connection. Gluing  these  two boundary  components  together  adds  a handle  to  the surface. So  this case  does not appear for genus zero.\\
 In the other case, the single boundary component  consists  of  two saddle  connections. The surface $S$ is obtained after gluing these two saddle connections, so $\gamma$ consists of a single saddle connection joining a singularity of order $k$ to a distinct singularity of order $k^{\prime}$. Note that $k$ and $k^{\prime}$ cannot be both equal to $-1$ otherwise there would be another saddle connection in the collection $\gamma$ (see remark \ref{generique}).
%

b) Two  ``$-$'' vertices  of valence one  joined by  a single edge. That means that $\gamma$ consists of a single closed saddle connection $\gamma_1$ which separates the surface in two parts. We get a unordered pair $\{S_1,S_2\}$ of one-holed spheres with  boundary singularities of angles $(a_1+1)\pi$ and $(a_2+1)\pi$ correspondingly. The saddle connection of the initial surface is adjacent to a singularity of order $a_1+a_2=k$. None of the $a_i$ is null otherwise the saddle connection would bound a cylinder, and there would exist a saddle connection \^homo\-logous to $\gamma_1$ on the other boundary component of this cylinder. 

 Now considering the interior singularities of positive order of $S_1$ and $S_2$ respectively, this defines a partition $A_1 \sqcup A_2$ of $\{k_1^{\alpha_1},\ldots,k_r^{\alpha_r}\}\backslash\{k\}$.
  Each $S_i$ also have $s_i$ poles, with $s_1+s_2=s$. If we decompose the boundary saddle connection of $S_i$ in two segments starting from the boundary singularity, and glue together these two segments,  we then get a  closed flat surface with $A_i\sqcup \{a_1-1,-1\}\sqcup \{-1^{s_i}\}$ for the order of the singularities. The Gauss-Bonnet theorem implies: \[\bigl(\sum_{a\in A_i} a\bigr)+a_1-2-s_i=-4.\]

c) Two  ``$-$'' vertices  of valence one  and a ``$\circ$'' vertex of valence $2$.
This case is analogous to the previous one.

d) A ``$-$'' vertex of valence  one, joined by an edge to a valence three  ``$\circ$'' vertex and  an edge joining the ``$\circ$'' vertex to itself.\\ 
The ``$-$'' vertex   represents a  one-holed sphere. It has  a  single boundary component  which is a  closed saddle  connection.  The cylinder has  two boundary components  of equal lengths. One  has two  saddle  connections of  length 1  (after  normalization)  the  other component has  a single saddle connection of  length $2$. So, the only possible  configuration is obtained by gluing the two saddle  connections of length  1 together (creating  a ``half-pillowcase'') and  gluing the other one with the boundary of the ``$-$'' component. The boundary singularity of the ``$-$'' component has an angle of $(k+2-1)\pi$ (equivalently, has order $k$) for some $k\in\{k_1,\ldots,k_r\}$. 

e) A valence four ``$\circ$'' vertex with two edges joining the vertex to itself. The cylinder has two boundary components,  each  of them  is composed  of  two  saddle connections. All the  saddle connections have the same length. If we glue a saddle connection with one of the other connected  component of  the boundary, we  get a flat torus, which has trivial holonomy and genus greater than zero. So, we have to glue  each  saddle connection  with  the other  saddle connection of its boundary component. That means that we get a (twisted) ``pillowcase'' and the surface belongs to $\mathcal{Q}(-1,-1,-1,-1)$.\\

In each of these first four cases, the surface necessary has a singularity of order at least one.  So, they cannot appear in $\mathcal{Q}(-1,-1,-1,-1)$, which means that the fifth case is the only possibility in that stratum.
\end{proof}

\section{Configurations for hyperelliptic connected components}

In this section, we describe the configurations of \^homo\-logous saddle
connections in a hyperelliptic connected component. We first reformulate
Lanneau's description of such components, see \cite{La1}.

\begin{NNths}[E. Lanneau]
The  hyperelliptic  connected components  are given by the following list:
\begin{enumerate}
\item  The subset of surfaces in $\mathcal{Q}(k_1,k_1,k_2,k_2)$, that are a double covering of a surface in  $\mathcal{Q}(k_1,k_2,-1^{s})$ ramified over $s$ poles. Here $k_1$ and $k_2$ are odd, $k_1\geq -1$ and $k_2\geq 1$, and $k_1+k_2- s=-4$.
\item  The subset of surfaces in $\mathcal{Q}(k_1,k_1,2k_2+2)$, that are a double covering of a surface in  $\mathcal{Q}(k_1,k_2,-1^{s})$ ramified over $s$ poles and over the singularity of order $k_2$.
Here $k_1$ is odd and  $k_2$ is even, $k_1\geq -1$ and $k_2\geq 0$, and
$k_1+k_2-s=-4$.
\item The subset of surfaces in $\mathcal{Q}(2k_1+2,2k_2+2)$, that are a double covering of a surface in  $\mathcal{Q}(k_1,k_2,-1^{s})$ ramified over all the singularities. 
Here $k_1$ and $k_2$ are even, $k_1\geq 0$ and $k_2\geq 0$, and $k_1+k_2-
s=-4$.
\end{enumerate}
\end{NNths}

\begin{table}[htbp]
   \begin{center}
 \input{confighyp2.pstex_t}
   \end{center}
   \caption{Configurations for $\mathcal{Q}^{hyp}(k_1, k_1 ,2k_2+2)$}
   \label{graphhyp3}
\end{table} 

\begin{table}[htbp]
   \begin{center}
 \input{confighyp1.pstex_t}
   \end{center}
   \caption{Configurations for $\mathcal{Q}^{hyp}(k_1,k_1,k_2,k_2)$}
   \label{graphhyp1}
\end{table}

\begin{table}[htbp]
   \begin{center}
       \input{confighyp3.pstex_t}
   \end{center}
   \caption{Configurations for $\mathcal{Q}^{hyp}(2k_1+2,2k_2+2)$}
   \label{graphhyp2}
\end{table}

Taking a double covering of the configurations arising on $\mathbb{CP}^1$, one can deduce configurations for hyperelliptic components. 
This leads to the following theorem:

\begin{ths} \label{confighyp}
In the notations of the classification theorem above, the admissible configurations of \^homo\-logous saddle connections for  hyperelliptic connected components are given by tables  \ref{graphhyp3}, \ref{graphhyp1}, \ref{graphhyp2} and \ref{graphhyp_mp}. No other configuration can appear.
\end{ths}

\begin{rem}
Integer parameters $k_1, k_2\geq -1$ in tables   \ref{graphhyp3}, \ref{graphhyp1}, \ref{graphhyp2}  are allowed to take values $-1$ and $0$ as soon as this does not contradict explicit restrictions. In table \ref{graphhyp_mp}, we list several additional configurations which appear only when at least one of $k_1,k_2$ is equal to zero.
\end{rem}
\begin{rem}
In the description of configurations for the hyperelliptic connected component $\mathcal{Q}^{hyp}(k_1,k_1,k_2,k_2)$ with $k_1=k_2$, the notation $k_i,k_i$ (\emph{resp.} $k_j,k_j$) still represents the orders of a pair of singularities that are \emph{interchanged} by the hyperelliptic involution. For example in a generic surface in the hyperelliptic component $\mathcal{Q}^{hyp}(k,k,k,k)$, for $k\geq 1$, the second line of table \ref{graphhyp1} means that,  between any pair of singularities that are interchanged by the hyperelliptic involution on $S$, there exists a saddle connection with no other saddle connections \^homo\-logous to it. But if $\gamma$ is a saddle connection between two singularities that are not interchanged by the involution $\tau$, then $\tau(\gamma)$ is  a saddle connection \^homo\-logous to $\gamma$ (see below), and which is different from $\gamma$.
\end{rem} 

\begin{proof}
Let $\mathcal{Q}^{hyp}$ be a hyperelliptic connected component as in the list of the previous theorem and $\mathcal{Q}=\mathcal{Q}(k_1,k_2,-1^s)$ the corresponding stratum on $\mathbb{CP}^1$. The projection $p:\tilde{S} \rightarrow \tilde{S}/\tau=S$, where $\tilde{S}\in \mathcal{Q}^{hyp}$ and $\tau$ is the corresponding hyperelliptic involution,  induces a covering  from $\mathcal{Q}^{hyp}$ to $\mathcal{Q}$. This is not necessarily a one-to-one map because there might be a choice of the ramification points on $\mathbb{CP}^1$. But if we fix the ramification points, there is a locally one-to-one correspondence.  

We recall to the reader that theorem of Masur and Zorich cited after definition \ref{def:homo} says that two saddle connections are \^homo\-logous if and only if the ratio of their length is constant under any small perturbation of the surface inside the ambient stratum. Therefore, two saddle connections in  $\tilde{S}\in \mathcal{Q}^{hyp}$ are \^homo\-logous if and only if the corresponding saddle connections in $S$ are \^homo\-logous. 
Hence  the image under $p$ of a maximal collection $\tilde{\gamma}$ of \^homo\-logous saddle connections on $\tilde{S}$ is a collection $\gamma$ of \^homo\-logous saddle connections on $S$. Note that $\gamma$ is not necessary maximal since the preimage of a pole by $p$ is a marked point on $\tilde{S}$ and we do not consider saddle connections starting from a marked point. However, we can deduce \emph{all} configurations for $\mathcal{Q}^{hyp}$ from the list of configurations for $\mathcal{Q}$.

We give details  for a few configurations, the other ones are similar and the proofs are left to the reader.

-\emph{First line of  table \ref{graphhyp3}}: the configuration for $\mathcal{Q}=\mathcal{Q}(k_1,k_2,-1^s)$ corresponds to a single saddle connection $\gamma$ on a surface $S$ that joins a singularity $P_1$ of degree $k_1$ to the distinct  singularity $P_2$ of degree $k_2$.  The double covering is ramified over $P_2$ but not over $P_1$. Therefore, the preimage  of $\gamma$ in $\tilde{S}$ is a pair $\{\tilde{\gamma}_1,\tilde{\gamma_2}\}$ of saddle connections of the same lengths that join each preimage of $P_1$ to the preimage of $P_2$.   The boundary  of compactification of $\tilde{S}\backslash \{\tilde{\gamma}_1,\tilde{\gamma_2}\}$ admits only one connected component that consists of four saddle connections. The angles of the boundary singularities corresponding to the preimages of $P_1$ are both $(k_1+2)\pi$, and the angles of the other boundary singularities are $(k_2+2)\pi$ since $\{\tilde{\gamma}_1,\tilde{\gamma_2}\}$ are interchanged by the hyperelliptic involution.

-\emph{Fourth line of table \ref{graphhyp3}}: the configuration for $\mathcal{Q}=\mathcal{Q}(k_1,k_2,-1^s)$ corresponds to a single closed saddle connection $\gamma$ on a flat surface $S$ that separates the surface into two parts $S_1$ and $S_2$. Each $S_i$ contains some ramification points, so the preimage of $\gamma$ separates $\tilde{S}$ into two parts $\tilde{S}_1$ and $\tilde{S}_2$ that are double covers of $S_1$ and $S_2$. One of the $\tilde{S}_i$ has an interior singularity of order $2k_2+2$, while the other one does not have interior singularities. The description from Masur and Zorich of possible graphs of connected components (see Figure \ref{fig:classification:of:graphs}) implies that $\tilde{S}_1$ and $\tilde{S}_2$ cannot have the same holonomy. Let $\tilde{S}_2$ be the component with trivial holonomy, and choose $\omega$  a square root of the quadratic differential that defines its flat structure. If $\tilde{S}_2$ has two  boundary components,  each consisting of a single saddle connection, then the corresponding boundary singularities must be of even order $a_2$. If $\tilde{S}_2$ has a single boundary component, then integrating $\omega$ along that boundary must give zero ($\omega$ is closed), which is only possible if the order $a_2$ of the boundary singularities are odd. Applying Lemma \ref{bord-revdouble} below, we see that $\tilde{S}_2$ does not have interior singularity. Hence, $\tilde{S}_1$ has an interior singularity of order $2k_2+2$.
The order of the boundary singularities of $\tilde{S}_1$ are both $a_1=k_1-a_2$, which is of parity opposite to the one of $a_2$.  Applying again Lemma \ref{bord-revdouble}, we get the number of boundary components of $\tilde{S}_1$.

-\emph{Last line of  table \ref{graphhyp3}}: the configuration for $\mathcal{Q}=\mathcal{Q}(k_1,k_2,-1^s)$ corresponds to a pair of saddle connections on a surface $S\in \mathcal{Q}$ that separate the surface into a cylinder $C$ and a one-holed sphere $S_1$. The double cover $\tilde{S}_1$ of $S_1$ is connected, and applying  Lemma \ref{bord-revdouble} we see that it has two boundary components. The double cover $\tilde{C}$ of the cylinder $C$ admits no ramification point. So \emph{a priori}, there are two possibilities: $\tilde{C}$ is either a cylinder  the same length of and a width twice bigger than the width of $C$, or it is a pair of copies of $C$. Here, the first possibility is not realizable otherwise the double covering $\tilde{S}\rightarrow S$ would be necessary ramified over $k_1$. Finally we get $\tilde{S}$ by gluing a boundary component of each cylinder to each boundary component of $\tilde{S}_1$, and gluing together the remaining boundary components of the cylinders. \\ 
Note that  the preimage of the saddle connection joining a pair of poles on $S$ is a \emph{regular} closed geodesic in $\tilde{S}$, and hence in our convention, we do not consider such a saddle connection in the collection of \^homo\-logous saddle connections on~$\tilde{S}$.

\bigskip

When at least one of $k_1$ or $k_2$ equals zero, there is a marked point on $\mathbb{CP}^1$ that is a ramification point of the double covering. Hence
we have to start from a configuration of saddle connections on $\mathbb{CP}^1$ that might have marked points as end points:
\begin{itemize}
\item If a maximal collection of \^homo\-logous saddle connection on $\mathbb{CP}^1$ does not intersect a marked point, then the collection has already been described in Theorem \ref{prop:config}, and hence, the corresponding configuration in $\mathcal{Q}^{hyp}$ is already presented in tables \ref{graphhyp3}  and \ref{graphhyp2}.
\item
If a non-closed saddle connection in a collection admits a marked point as end point, then this saddle connection is simple since we can move freely that marked point. Hence the corresponding configuration in $\mathcal{Q}^{hyp}$ is already written in tables \ref{graphhyp3}  and \ref{graphhyp2}.
\item 
If a closed saddle connection admits a marked point as end point, then it is a closed geodesic. This corresponds to a new configuration on $\mathbb{CP}^1$ and the corresponding configuration in $\mathcal{Q}^{hyp}$ is described in table \ref{graphhyp_mp}. The proof is analogous to the other cases.
\end{itemize}

\begin{table}[htbp]
   \begin{center}
       \input{confighyp_marked_point.pstex_t}
   \end{center}
   \caption{Additional configurations which appears when at least one of $k_1$ or $k_2$ equals $0$.}
   \label{graphhyp_mp}
\end{table}

This concludes the proof of Theorem \ref{confighyp}.

\end{proof}

\begin{lem}\label{bord-revdouble}
Let $S_{i}$ be a flat surface whose boundary consists of a single
closed saddle connection and let $a>0$ be the order of the corresponding
boundary singularity. Let $\tilde{S}_i$ be a connected ramified double cover of
the interior of $S_{i}$ and  let $(\tilde{k}_1,\ldots,\tilde{k}_l)$ be interior singularities. The sum $\sum_i \tilde{k}_i$ is even and:
\begin{itemize}
\item
If $\frac{\sum_i \tilde{k}_i}{2}+a$ is even, then the compactification of
$\tilde{S}_i$ has two boundary components, each of them consists of a single
saddle connection, with corresponding boundary singularity of order $a$.
\item If $\frac{\sum_i \tilde{k}_i}{2}+a$ is odd, then the
compactification of $\tilde{S}_i$ has a single boundary component which consists of a pair of saddle connections of equal lengths,  with corresponding boundary singularities of order $a$.
\end{itemize}
\end{lem}

\begin{proof}
By construction, the boundary of the compactification of $\tilde{S}_i$
necessary consists of two saddle connections of equal lengths. It has one or
two connected components.

Now we claim that 
\[\sum_i \tilde{k}_i + 2a \equiv 2r \mod 4 \]
 where $r$ is the
number of connected components of the boundary of~$\tilde{S}_i$. This equality
(that already appears in  \cite{MZ}) clearly implies the lemma. To prove the claim, we consider as in Lemma \ref{genre_plus} the surface $\tilde{X}$ of genus $g_{\tilde{X}}$ obtained by gluing $\tilde{S}_i$ and a copy of itself with opposite orientation along their boundaries. The orders of the singularities of  $\tilde{X}$ are 
 $\{\tilde{k}_1,\ldots,\tilde{k}_l,\tilde{k}_1,\ldots,\tilde{k}_l,2a,2a\}$, so we get 
\[4g_{\tilde{X}}-4=2\sum_i \tilde{k}_i +4a = 4(2\tilde{g}_i +r-1)-4 \]
 and  therefore \[\sum_i \tilde{k}_i +2a = 4g_i-4 +2r \equiv 2r \mod 4.\] 
\end{proof}

Given a concrete flat surface, we do not necessary see at once whether it belongs or not to a hyperelliptic connected component. Indeed, there exists hyperelliptic flat surfaces that are not in a hyperelliptic connected component.
As a direct corollary of Theorem \ref{confighyp}, we have the following quick test.

\begin{cor}
Let $S$ be a flat surface with non-trivial holonomy and let $\gamma$ be a collection of \^homo\-logous saddle connections on $S$.
If one of the following property holds, then the surface $S$ does not belong to a hyperelliptic connected component.
\begin{itemize}
\item  $S\backslash \gamma$ admits three connected components and neither of them is a cylinder.
\item $S\backslash \gamma$ admits four connected components or more.
\end{itemize}
\end{cor}

\section{Configurations for non-hyperelliptic connected components}

Following \cite{MZ}, given a fixed stratum, one can get a list of all realizable configurations of \^homo\-logous saddle connections. Nevertheless it is not clear which configuration realizes in which component. In the previous section we have described configurations for hyperelliptic components.

In the section we show that \emph{any} configuration realizable for a stratum is  realizable in its non-hyperelliptic connected component, provided the genus $g$ is sufficiently large.

We will use the following theorem which is a reformulation of the theorem of Kontsevich-Zorich and the theorem of Lanneau cited in section \ref{bas_def}. 

\begin{NNths}[M. Kontsevich, A. Zorich; E. Lanneau]
The following strata consists entirely of hyperelliptic surfaces and are connected.
\begin{itemize}
\item $\mathcal{H}(0)$, $\mathcal{H}(0,0)$, $\mathcal{H}(1,1)$ and $\mathcal{H}(2)$ in the moduli spaces of Abelian differentials.
\item $\mathcal{Q}(-1,-1,1,1)$, $\mathcal{Q}(-1,-1,2)$, $\mathcal{Q}(1,1,1,1)$, $\mathcal{Q}(1,1,2)$ and $\mathcal{Q}(2,2)$ in the moduli spaces of quadratic differentials.
\end{itemize}
Any other stratum that contains a hyperelliptic connected component admit at least one other connected component that contains a subset of full measure of flat surfaces that do not admit any isometric involution.
\end{NNths}

\begin{lem}\label{simple_dans_nonhyp}
Let $\mathcal{Q}$ be a non-connected stratum that  contains a hyperelliptic connected component. If the set of order of singularities defining $\mathcal{Q}$ contains $\{k,k\}$, for some $k\geq 1$, then there exists a non-hyperelliptic flat surface in $\mathcal{Q}$ having a simple saddle connection joining two different singularities of the same order $k$.
\end{lem}

Here we call a saddle connection ``simple'' when there are no other saddle connections \^homo\-logous to it.

\begin{proof}
According to Masur and Smillie \cite{MS}, any stratum is nonempty except the following four exceptions: $\mathcal{Q}(\emptyset)$, $\mathcal{Q}(1,-1)$, $\mathcal{Q}(3,1)$ and $\mathcal{Q}(4)$. 

According to Masur and Zorich \cite{MZ} (see also \cite{EMZ}), if $S\in \mathcal{Q}(k_1+k_2,k_3,\ldots,k_r)$, then there is a continuous path $(S_t)_{t\in [0,1]}$ in the moduli space of quadratic differentials, such that $S_0=S$ and $S_t$ is in $\mathcal{Q}(k_1,k_2,k_3,\ldots,k_r)$ for $t>0$, and such that the smallest saddle connection on $S_t$, for $t>0$ is simple and joins a singularity of order $k_1$ to a singularity of order $k_2$.  We say that we ``break up'' the singularity of order $k_1+k_2$ into two singularities of order $k_1$ and $k_2$.

\begin{figure}[htb]
\begin{center}
\input{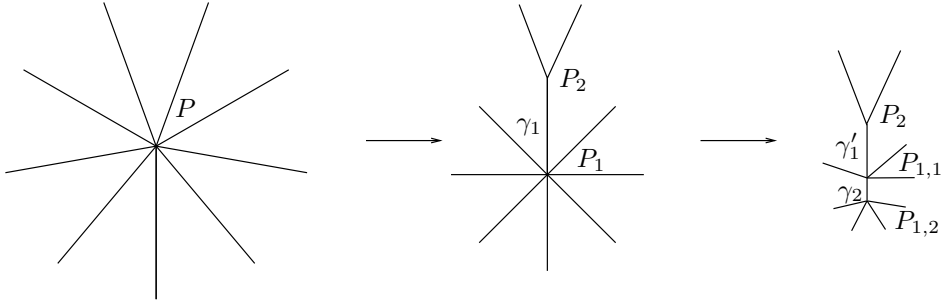}
\caption{Construction of  a simple saddle connection in a non-hyperelliptic surface}
\label{simple_non_hyp}
\end{center}
\end{figure}

We first consider the stratum $\mathcal{Q}=\mathcal{Q}(k_1,k_1,k_2,k_2)$.
By assumption, $\mathcal{Q}$ is non-connected, so, either the genus is greater than $3$, or $k_1=3$ and $k_2=-1$. Hence the stratum $\mathcal{Q}(2k_1+k_2,k_2)$ is nonempty. Now, we start from a surface $S_0$ in that stratum, and break up the singularity $P$ of order $2k_1+k_2$ into two singularities $P_1$ and $P_2$ of orders $2k_1$ and $k_2$ respectively (see Figure \ref{simple_non_hyp}). We get a surface $S_1$ with a short vertical saddle connection $\gamma_1$ between  $P_1$ and $P_2$. Since the ``singularity breaking up'' procedure is continuous, there are no other short saddle connections on $S_1$. Then, we break up the singularity $P_1$ of order $2k_1$ into a pair of singularities $P_{1,1}$ and $P_{1,2}$ of orders $k_1$. We get by construction a surface $S_2$ in the stratum $\mathcal{Q}$ with a simple saddle connection $\gamma_2$ between $P_{1,1}$ and $P_{1,2}$, and of length very small compared to the length of $\gamma_1$. 

The fact that the ``singularity breaking up'' procedure is continuous implies that there persists a saddle connection $\gamma_1^{\prime}$ between $P_2$ and one of the $P_{1,i}$ (see Figure \ref{simple_non_hyp}). By construction, we can assume there is no other saddle connection of length $\kappa l(\gamma_1^{\prime})$, where $l(\gamma_1^{\prime})$ denotes the length of $\gamma_1^{\prime}$ and $\kappa\in \{\frac{1}{2},1,2\}$ . Hence, $\gamma_1^{\prime}$ is simple by theorem of Masur and Zorich cited after definition \ref{def:homo}. According to Theorem \ref{confighyp}, this cannot exist in the hyperelliptic connected component since the corresponding configuration is not present in table \ref{graphhyp1}. Thus $S_2$ belongs to the non-hyperelliptic connected component and we can assume, after a slight perturbation, than $S_2$ is not hyperelliptic. Since by construction, the saddle connection $\gamma_2$ is simple and joins two singularities of order $k=k_1\geq 1$, the lemma is proven for the stratum $\mathcal{Q}(k_1,k_1,k_2,k_2)$.

The proofs for $\mathcal{Q}(k_1,k_1,2k_2+2)$ and for $\mathcal{Q}(2k_1+2,2k_2+2)$ are analogous: note that these case do not occur for the genera $1$ or $2$, because all corresponding strata are connected. Therefore the genus is greater than or equal to $3$ and the stratum  $\mathcal{Q}(2k_1+2k_2+2)$ is nonempty.
\end{proof}

\begin{ths}\label{config:non:hyp}
Let $\mathcal{Q}$ be a stratum of meromorphic quadratic differentials with at most simple poles on a Riemann surface of genus $g\geq 5$. If $\mathcal{Q}$ admits a hyperelliptic connected component, then $\mathcal{Q}$ is non-connected and any configuration for $\mathcal{Q}$ is realized for a surface in the non-hyperelliptic connected component of $\mathcal{Q}$.
\end{ths}

\begin{proof}
The fact that $\mathcal{Q}$ is non-connected follows directly from the Theorem of Lanneau.
Let $S$ be a flat surface in the hyperelliptic component of $\mathcal{Q}$ and let $\gamma$ be a maximal collection of \^homo\-logous saddle connections.  The hyperelliptic involution $\tau$ maps $\gamma$ to itself and hence, induces an involution $\tau^*$ on the set of connected components of $S\backslash \gamma$. 
We claim that $\tau^*$ does not interchange two connected components $S_i, S_j$ of $S\backslash \gamma$, for otherwise we can continuously deform $S_i$ outside a neighborhood of its boundary and 
 reconstruct a new flat surface $S^{\prime}$ in the same connected component.
 By construction, such surface $S^{\prime}$ is not any more hyperelliptic. Therefore, if $S$ is in a hyperelliptic component, then $\tau$ must induce an isometric and orientation preserving involution on each connected component of $S\backslash \gamma$.

Using the formula for the genus of a compound surface proved in the appendix and the list of configurations for hyperelliptic connected components given in the previous section, we derive the following fact: if $S$ has genus $g\geq 5$ and if $\gamma$ is a maximal collection of \^homo\-logous saddle connections, then at least one of the following three propositions is true.
\begin{itemize}
\item[a)]  $S\backslash \gamma$ admits a connected component $S_{0}$ of genus $g_{0} \geq 3$, that has a single boundary component and 
whose corresponding vertex in the graph $\Gamma(S,\gamma)$ is of valence $2$.

\item[b)]  $S\backslash \gamma$ admits a connected component $S_{0}$ of genus $g_{0} \geq 2$, that has exactly two boundary components and 
whose corresponding vertex in the graph $\Gamma(S,\gamma)$ is of valence $2$.

\item[c)] $S\backslash \gamma$ is connected and the corresponding vertex in the graph $\Gamma(S,\gamma)$ is of valence $4$. 
\end{itemize}

The proof follows from Lemmas \ref{lem_b}, \ref{lem_a}, \ref{lem_c} to situations a), b), c) correspondingly.
\end{proof}

\begin{lem}  \label{lem_b}
Let $S$ be a flat surface in a hyperelliptic connected component and let $\gamma$ be a maximal collection of \^homo\-logous saddle connections. We assume that 
$S\backslash \gamma$ admits a connected component $S_{0}$ of genus $g_{0} \geq 3$, whose corresponding vertex in the graph $\Gamma(S,\gamma)$ is of valence $2$, and such that $S_0$ has a single boundary component. 

Then there exists $(S^{\prime}, \gamma^{\prime})$ that has the same configuration as $(S,\gamma)$, with $S^{\prime}$ in the complementary component of the same stratum.
\end{lem}

\begin{proof}
The boundary components of $S_0$  consists of two saddle connections of the same length and the corresponding boundary singularities have the same  orders $k\geq 1$. Identifying together these two boundary saddle connections, we get a hyperelliptic surface $\overline{S}_0$. If we continuously deform this surface, it keeps being hyperelliptic since we can perform the reverse surgery and get a continous deformation of $S$.
 Hence, $\overline{S}_0$ belongs to a hyperelliptic component, and the hyperelliptic involution interchange two singularities of order $k-1$.

The genus of $\overline{S}_0$ is greater than $3$, so the corresponding stratum admits an other connected component. Now we start from a closed flat surface $X$ in this other connected component.   According to Lemma \ref{simple_dans_nonhyp},  we  can choose $X$ such that it admits a simple saddle connection between the two singularities of order $k-1$.
Now we cut $X$ along that saddle connection and we get a surface $S_1$ that have, after rescaling, the same boundary as $S_0$. By construction, $S_1$ admits no interior saddle connection \^homo\-logous to one of its boundary  saddle connections. So, we can reconstruct a pair $(S^{\prime}, \gamma^{\prime})$ such that $\gamma^{\prime}$ has the same configuration as $\gamma$ in $S$.

The surface $S_1$ admits a nontrivial isometric involution if and only if $X$ shares this property. So,  we can choose $X$ in such a way it admits no nontrivial isometric involution, and therefore the surface $S^{\prime}$ is non-hyperelliptic.

This argument also works when $\overline{S}_0$ is in the stratum $\mathcal{Q}(3,3,-1,-1)$ (here $g_0=2$ and $k=4$). In any other case for $g_0\leq 2$, it is not possible to replace $S_0$ by a surface $S_1$ with no involution.
\end{proof}

\begin{lem} \label{lem_a}
Let $S$ be a flat surface in a hyperelliptic connected component and let $\gamma$ be a maximal collection of \^homo\-logous saddle connections. We assume that $S\backslash \gamma$ admits a connected component $S_{0}$  of  genus $g_{0}\geq 2$, that has two boundary components, and  whose corresponding vertex in the graph $\Gamma(S,\gamma)$ is of valence $2$. 

Then there exists $(S^{\prime}, \gamma^{\prime})$ that has the same configuration as $(S,\gamma)$, with $S^{\prime}$ in the complementary component of the same stratum.
\end{lem}

\begin{proof}
Each boundary component of $S_0$  consists of one saddle connection and the corresponding boundary singularities have the same  orders $k\geq 1$. Now we start from a closed flat surface $X$ with the same holonomy as $S_{0}$ and whose singularities consists of the interior singularities of $S_{0}$ and two singularities $P_{1}$ and $P_{2}$ of order $k-2$. We can always choose $X$ such that it admits a saddle connection $\eta$ between $P_{1}$ and~$P_{2}$. 

 Now we construct a pair of holes by removing a parallelogram as in Figure $\ref{a_case}$  and gluing together the two long  sides. 
 Note that the holes can be chosen arbitrary small, and therefore, the resulting surface with boundary does not have any interior saddle connection \^homo\-logous to one of its boundary components. We denote by $S_1$ this surface, and up to rescaling, we can assume that $S_0$ and $S_1$ have isometric boundaries.  Hence replacing $S_0$ by $S_1$ in the decomposition of $S$, we get a new pair $(S^{\prime},\gamma^{\prime})$ that have the same configuration as $(S,\gamma)$. 

\begin{figure}[htbp]
   \begin{center}
       \input{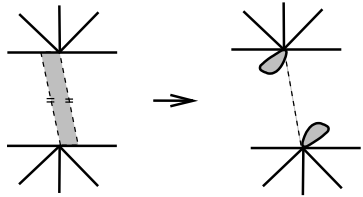}
   \end{center}
   \caption{Construction of a pair of holes}
   \label{a_case}
\end{figure}

We now assume that $S_1$ admits a nontrivial isometric (orientation preserving) involution $\tau$. Then this involution interchanges the two boundary components of the surface. It is easy to check that  we can perform the reverse surgery as the one described previously and we get a closed surface that admits a nontrivial involution. Hence if $X$ belongs to a stratum that does not consist entirely of hyperelliptic flat surfaces, then we can choose $X$ such that $S^{\prime}$ is not in a hyperelliptic connected component.

The hypothesis on the genus, the  theorem of Kontsevich-Zorich and the theorem of Lanneau imply that this arguments works except when $X$ belongs to $\mathcal{H}(1,1)$, $\mathcal{Q}(2,1,1)$, $\mathcal{Q}(1,1,1,1)$, or $\mathcal{Q}(2,2)$. 

We remark that if $X\in \mathcal{Q}(2,2)$, then $S_0$ must have nontrivial linear holonomy and no interior singularities. According to the list of configurations for hyperelliptic connected components given in section $3$, this cannot happen.

We exhibit in Figure \ref{b_case} three explicit surfaces with boundary that corresponds to the three cases left. We represent these three surfaces as having a one-cylinder decomposition and by describing the identifications on the boundary of that cylinder. The length parameters can be chosen freely under the obvious condition that the sum of the lengths corresponding to the top of the cylinder must be equal to the sum of the lengths corresponding to the bottom of the cylinder. Bold lines represents the boundary of the flat surface. Now we remark that a nontrivial isometric involution must preserve the interior of the cylinder, and must exchange the boundary components. This induces some supplementary relations on the length parameters. Therefore, we can choose them such that there is no nontrivial isometric involution.

\begin{figure}[htbp]
   \begin{center}
       \input{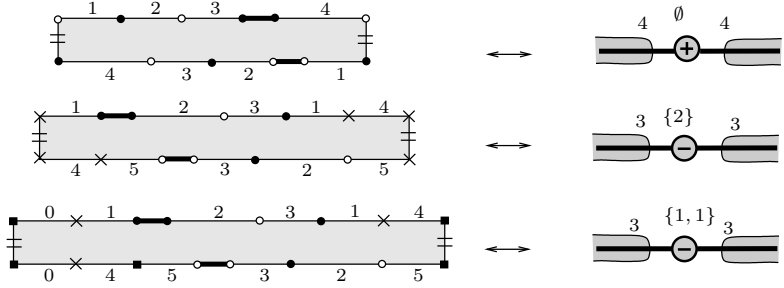}
   \end{center}
   \caption{Surfaces with two boundary components and no involution in low genus.}
   \label{b_case}
\end{figure}

\end{proof}

\begin{lem} \label{lem_c}
Let $S$ be a flat surface of genus $g\geq 3$ with nontrivial linear holonomy that belongs to a hyperelliptic connected component and let $\gamma=\{\gamma_1,\gamma_2\}$ be a maximal collection of \^homo\-logous saddle connections on $S$. If $S\backslash \gamma$ is connected, 
then there exists $(S^{\prime}, \gamma^{\prime})$ that has the same configuration as $(S,\gamma)$, with $S^{\prime}$ in the complementary component of the same stratum.
\end{lem}

\begin{figure}[htb]
   \begin{center}
       \input{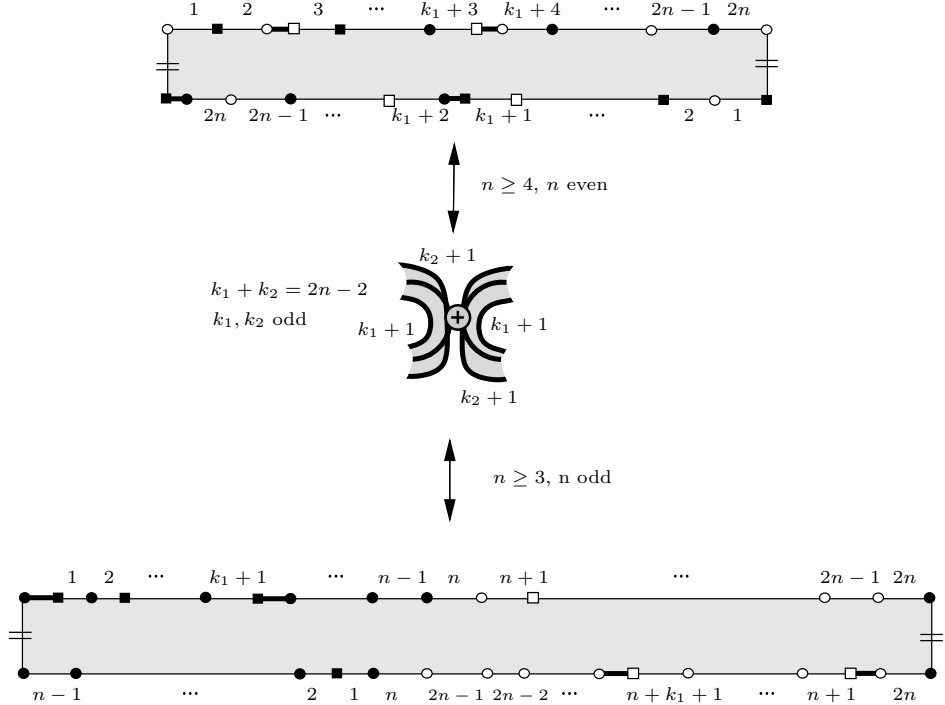}
   \end{center}
   \caption{Valence four component with no involution}
   \label{c_case}
\end{figure}

\begin{proof}
Since $S\backslash\gamma$ is connected,  the graph $\Gamma(S,\gamma)$ contains a single vertex, and it has valence four. According to Theorem \ref{confighyp}, two different cases appear:

a) The surface $\overline{S\backslash \gamma}$ has one boundary component. In this case we start from a surface in $\mathcal{H}(k_1+k_2+1)$ and perform a local surgery in a neighborhood of the singularity, as described in Figure \ref{local4} (see also \cite{MZ}, section $5$). We get  a surface and a pair of small saddle connections of length $\delta$ that have the same configuration as $\gamma$. The stratum $\mathcal{H}(k_1+k_2+1)$ admits non-hyperelliptic components and the  same argument as in the previous lemmas works: if we start from a generic surface in a non-hyperelliptic component, then the resulting surface after surgery does not have any nontrivial involution.

\begin{figure}[htbp]
   \begin{center}
       \input{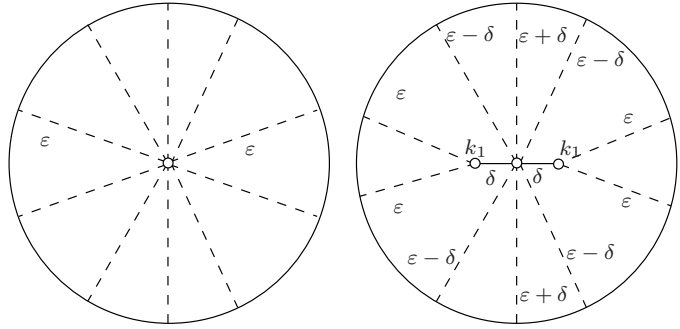}
   \end{center}
   \caption{Breaking up a zero in three ones}
   \label{local4}
\end{figure}

b) The surface $\overline{S\backslash \gamma}$ has two boundary components, each of them consists of a pair of saddle connections with boundary singularities of order $k_1+1$ and $k_2+1$. We construct explicit surfaces with the same configuration as $\gamma$, but that have no nontrivial involution. Let $2n=k_1+k_2+2$ and 
we start from a surface $S_0$ of genus $n$ in $\mathcal{H}(n-1,n-1)$, that have a one-cylinder decomposition and such as identification on the boundary of that cylinder is given by the permutation
\[
\left(\begin{array}{cccc}1 & 2 & \ldots & 2n \\ 2n & 2n-1 & \ldots  & 1\end{array}\right)\] 
when $n$ is even, and otherwise by the permutation 
\[
\left(\begin{array}{cccccccccc}
1 & 2 & \ldots &n-1 & n & n+1 &n+2 & \ldots & 2n-1 & 2n  \\
n-1&n-2 & \ldots & 1 & n & 2n-1 & 2n-2 & \ldots & n+1 & 2n
\end{array}\right)\]

We assume that $k_1$ and $k_2$ are odd and we perform a surgery on $S_0$ to get a surface $S_1$ with boundary as pictured on Figure \ref{c_case}. The surface $S_1$ admits two boundary components that consist of two saddle connections each  and which are represented by the bold segments. Each symbol \includegraphics{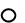}, \includegraphics{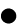}, \includegraphics{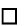}, \includegraphics{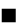} represents a different boundary singularity. It is easy to check that the boundary angles corresponding to  \includegraphics{carre_noir} and \includegraphics{carre_blanc} are both $(k_1+2)\pi$ and that the angles corresponding to \includegraphics{rond_noir} and \includegraphics{rond_blanc} are $(k_2+2)\pi$. Hence after suitable identifications of the boundary of $S_1$, we get a surface $S^{\prime}$ and a pair of \^homo\-logous saddle connections $\gamma^{\prime}$ that have the same configuration as $(S,\gamma)$.
However, $S^{\prime}$ does not admit any nontrivial involution if the length parameters are chosen generically.
Note that this construction does not work when $n=2$, but according to section $3$, and since $k_1$ and $k_2$ are odd, we have $n=g$, which is  greater than or equal to $3$ by assumption.

The case $k_1$ and $k_2$ even is analogous and left to the reader (note that in this case, $g=n+1$, and the construction works also for $n=2$).
\end{proof}

\appendix
\section*{Appendix.  Computation of the genus  in terms of a configuration}

Here we improve Lemma \ref{genre} and give the relation between the genus of a surface and the genera of the connected components of $S\backslash \gamma$, where $\gamma$ is a collection of \^homo\-logous saddle connections.

We first remark that this relation depends not only on the graph of connected components, but also on the permutation on each of its vertices (\emph{i.e.} on the ribbon graph). Indeed, let us consider a pair of \^homo\-logous saddle connections that decompose the surface into two  connected components $S_1$ and $S_2$. Then either both $S_1$ and $S_2$ have only one boundary component, or at least one of them  has two boundary components.  In the first case, $S$ is the connected sum of $\tilde{S}_1$ and $\tilde{S}_2$, so $g=g_1+g_2$, while in the second case, one has $g=g_1+g_2+1$.

\begin{APdefs}
Let $(S,\gamma)$ be a flat surface with a collection of \^homo\-logous saddle connections. The \emph{pure} ribbon graph associated to $(S,\gamma)$ is the 2-dimensional topological manifold obtained  from the ribbon graph by forgetting the graph $\Gamma(S,\gamma)$, as in Figure \ref{purerg}.

\begin{figure}[ht]
\begin{center}
\input{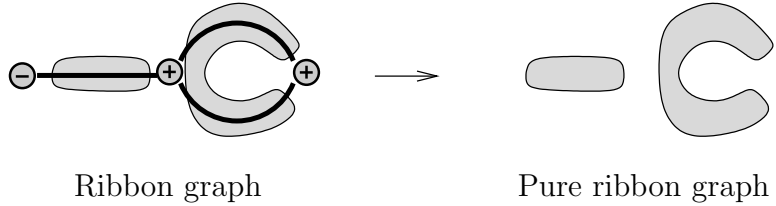}
\caption{Pure ribbon graph}
\label{purerg}
\end{center}
\end{figure}

\end{APdefs}

\begin{figure}[htbp]
\begin{center}
\includegraphics{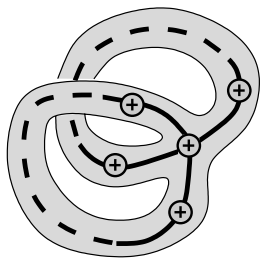}
\caption{Example of a ribbon graph that does not embed into $\mathbb{R}^2$.}
\label{do_not_embed}
\end{center}
\end{figure}

\begin{APprop}\label{genre_exact}
Let $\chi_1$ be the Euler characteristic of $\Gamma(S,\gamma)$ and let $\chi_2$ be the Euler characteristic of the pure ribbon graph associated to the configuration.
\begin{itemize}
\item If the pure ribbon graph has only one connected component and  does not embed into the plane (see Figure \ref{do_not_embed}), then 
\[
g=\bigl(\sum_i g_i \bigr)+1\] 
\item In any other case, 
\[
g=\bigl(\sum_i g_i \bigr)+(\chi_2-n)-(\chi_1-1)\]
\end{itemize}
\end{APprop}

\begin{APrem}
Simply connected components of the pure ribbon graph do not contribute to the term $(n-\chi_2)$, since the Euler characteristic of a disc is $1$.

Note also that  in the first case, we have $\chi_1=-1$ and $\chi_2=-1$, and therefore $\bigl(\sum_i g_i \bigr)+1\neq \bigl(\sum_i g_i \bigr)+(\chi_2-n)-(\chi_1-1)$. 
\end{APrem}

\begin{figure}[htbp]
\begin{center}
\input{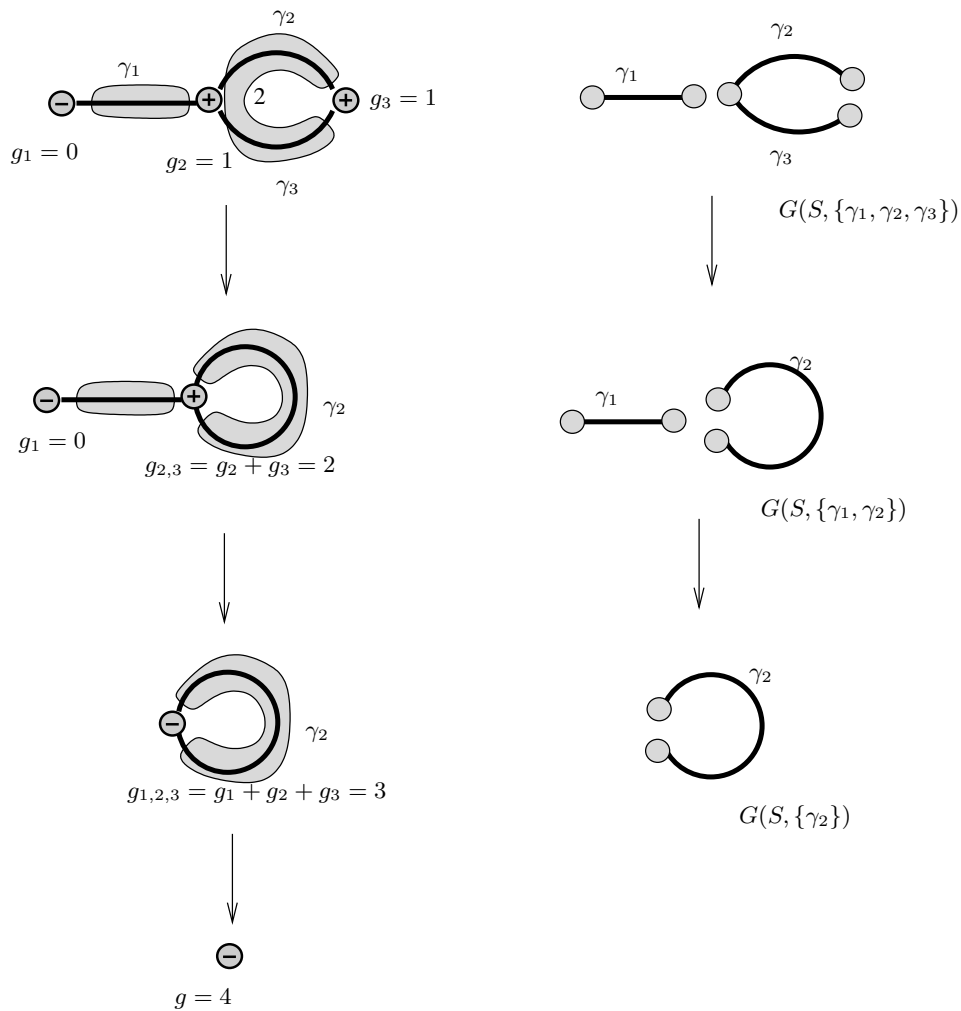}
\caption{Removing successively some elements of a collection ($\gamma_1,\gamma_2,\gamma_3)$.}
\label{exemple_genre_exact}
\end{center}
\end{figure}

\begin{proof}

Here we do not assume that the collection $\gamma$ is necessary maximal.
When $\Gamma(S,\gamma)$ has a single vertex, then we prove the proposition using direct computation and the description of the boundary components corresponding to each possible ribbon graph. We refer to \cite{MZ} for this description.
Then our goal is to reduce ourselves to that case by removing successively from the collection $\gamma=\{\gamma_1, \ldots,\gamma_k\}$ some $\gamma_i$ whose corresponding edges joins a vertex to a distinct one.

We define a new graph $G(S,\gamma)$, which is a deformation retract of the pure ribbon graph: the vertices of $G(S,\gamma)$ are the boundary components of each $S_i$, while the edges correspond to the saddle connections in $\gamma$ (see Figure \ref{exemple_genre_exact}). For each vertex, there is a cyclic order on the set of edges adjacent to the vertex consistent with the orientation of the plane. If the initial pure ribbon graph does not embed into the plane, then it is also the case for $G(S,\gamma)$. By construction, the Euler characteristic of $G(S,\gamma)$ is the same as the pure ribbon graph associated to $(S,\gamma)$, and is easier to compute.

Let $\Gamma(S,\gamma)$ contains at least two vertices. Choose a saddle connection representing an edge joining two distinct vertices of $\Gamma(S,\gamma)$, and up to renumeration, we can assume that this saddle connection is $\gamma_1$. Let us study the resulting configuration of $\gamma^{\prime}=\gamma\backslash \{\gamma_1\}$. The saddle connection $\gamma_1$  is on the boundary of two surfaces $S_{1}$ and $S_{2}$. Then the connected components of $S\backslash\gamma^{\prime}$ are the same as the connected component of $S\backslash\gamma$ except that the surfaces $S_{1}$ and $S_{2}$ are now glued along $\gamma_1$, and hence define a single surface $S_{1,2}$. The genus of $S_{1,2}$ (after gluing disks on its boundary) is $g_{1}+g_{2}$.

The graph $G(S,\gamma^{\prime})$ is obtained from $G(S,\gamma)$ by shrinking an edge that joins two different vertices, so these two graphs have the same Euler characteristic $\chi_1$.

Furthermore, if $\gamma_1$ was in a boundary component of $S_1$ (resp. $S_2$) defined by the ordered collection $(\gamma_1, \gamma_{i_1},\dots,\gamma_{i_s})$ 
(resp. $(\gamma_1,\gamma_{j_1},\dots,\gamma_{j_t})$). Then the  cyclic order in the corresponding boundary component of $S_{1,2}$ is defined by $(\gamma_{i_1},\dots,\gamma_{i_s}, \gamma_{j_1},\dots,\gamma_{j_t})$.
Therefore
$G(S,\gamma^{\prime})$ is obtained from $G(S,\gamma)$ by shrinking the edge corresponding to $\gamma_1$ and removing an isolated vertex that might appear (see Figure \ref{exemple_genre_exact}). It is clear that the difference $(\chi_2-n)$ between the Euler characteristic of $G(S,\gamma)$ and its number of connected component is constant under this procedure.  One can also remark that if $G(S,\gamma)$ is connected and does not embed into the plane (case $1$ of the proposition), then this is also true for $G(S,\gamma^{\prime})$.

Forgetting successively these $\gamma_i$ will lead to the case when $\Gamma(S,\gamma)$ has a single vertex. At each steps of the removing procedure, the numbers $\chi_1$ and $\chi_2-n$ do not change, and the sum of the genera associated to the vertices does not change either. This concludes the proof.

\end{proof}

\end{document}